\documentclass[11pt]{amsart}

\newcommand{\vs}{\vspace{.15in}}

\def\de{{\partial}}

\def\C{{\mathbf C}}
\def\P{{\mathbf P}}

\def\G{\Gamma}

\def\bX{\bar {X}}

\def\cS{{\mathcal S}}
\def\cP{{\mathcal P}}

\def\Z{{\bf Z}}
\def\Q{{\bf Q}}

\renewcommand{\|}{/\negmedspace/}

\newtheorem{theorem}{Theorem}[section]
\newtheorem{lemma}[theorem]{Lemma}
\newtheorem{proposition}[theorem]{Proposition}
\newtheorem{corollary}[theorem]{Corollary}

\theoremstyle{definition}
\newtheorem*{remark}{Remark}

\title{$\C_+$--actions on contractible threefolds}
\thanks{The first author was partially supported by NSA Grant
MDA904-00-1-0016. The second author was partially supported
by NSF Grant 0196523.}
\author[Shulim Kaliman]{Shulim Kaliman}
\author[Nikolai Saveliev]{Nikolai Saveliev}
\address{Department of Mathematics\newline\indent
University of Miami \newline\indent PO Box 249085
\newline\indent Coral Gables, FL 33124}
\email{\rm{kaliman@math.miami.edu}}
\email{\rm{saveliev@math.miami.edu}}

\begin{document}

\begin{abstract}  Let $X$ be a smooth contractible affine
algebraic threefold with a nontrivial algebraic $\C_+$-action on
it. We show that $X$ is rational and the algebraic quotient
$X\|\C_+$ is a smooth contractible surface $S$ which is isomorphic
to $\C^2$ in the case when $X$ admits a dominant morphism from a
threefold of form $C \times \C^2$. Furthermore, if the action is
free then $X$ is isomorphic to $S \times \C$ and the action is
induced by translation on the second factor. In particular, we
have the following criterion: if a smooth contractible affine algebraic
threefold $X$ with a free algebraic $\C_+$-action admits a
dominant morphism from $C\times \C^2$ then $X$ is isomorphic to
$\C^3$.
\end{abstract}

\maketitle

\section{Introduction} The aim of this paper is to generalize the
theorem of Miyanishi [Miy80] which says that, for any non-trivial
algebraic $\C_+$--action on $\C^3$, the algebraic quotient
$\C^3{\|}\C_+$ is isomorphic to $\C^2$. Our main result is that,
for a nontrivial algebraic $\C_+$-action on a smooth contractible
affine algebraic threefold $X$, the algebraic quotient $X\|\C_+$
is isomorphic to a smooth contractible affine surface $S$. As
all such surfaces are rational [GuSh89], we deduce that $X$ is
rational as well. Furthermore, if the action is free, we conclude
that $X$ is isomorphic to $S \times \C$ and the action is induced
by translation on the second factor, by virtue of [Ka03] where
this result was proved under the additional assumption that $S$ is
smooth. Another consequence of our main result is that, when $X$
admits a dominant morphism from a threefold of form $C \times
\C^2$, the quotient $S$ is isomorphic to $\C^2$. We also give an
independent proof of the latter fact which, unlike our main
result, does not use the difficult theorem of Fintushel and Stern
[FiSt90] about the absence of simply connected homology cobordisms
between certain homology spheres. In fact, the rationality of $X$
can also be proved without this theorem; however, this would
require another difficult theorem that all logarithmic
$\Q$-homology planes are rational [PrSh97, GuPrSh97, GuPr99]. In
conclusion we derive the following criterion : if there is a free
algebraic $\C_+$--action on a smooth contractible affine algebraic
threefold $X$ which admits a dominant morphism from $C\times \C^2$
then $X$ is isomorphic to $\C^3$.

\section{The main result}

Let $\rho: X \to S$ be the quotient morphism of a nontrivial
algebraic $\C_+$--action on a smooth contractible affine algebraic
threefold $X$. By Fujita's result, $X$ is factorial,
see e.g.[Ka94]. Some other properties of $\rho: X\to S$ proved in
[Ka03, Lemma 2.1, Proposition 3.2, Remark 3.3] are summarized in
the following lemma.

\begin{lemma}\label{lemma2.1} $\empty$

\begin{enumerate}
\item[(1)] The surface $S$ is affine and factorial and
$\rho^{-1}(s)$ is a nonempty curve for every $s\in S$.
\item[(2)] There is a curve $\G$ in $S$ so that $\breve{S} =
S\setminus \G$ is smooth and $\rho^{-1}( \breve{S})$ is
naturally isomorphic to $\breve{S } \times \C$ so that the
projection onto the first factor corresponds to $\rho$.
\end{enumerate}
\end{lemma}

\medskip

\begin{lemma}\label{lemma2.2}
In the above notation, let $S^*$ be the smooth part of the
quotient $S = X\|\C_+$. Then the groups $\pi_1 (S^*)$ and $H_2
(S^*)$ are trivial.
\end{lemma}

\begin{proof}
The set $F$ of singular points of $S$ is finite as $S$ is factorial.
According to Lemma \ref{lemma2.1}, $L = \rho^{-1}(F)$ is a curve,
hence $\pi_1 (X\setminus L) =\pi_2 (X\setminus L) = 0$.

Let $\gamma$ be a loop in $S^* = S\setminus F$. After a small
homotopy if necessary, we may assume that $\gamma \subset
\breve{S}$ where $\breve{S}$ is as in Lemma \ref{lemma2.1}. Since
$\rho^{-1}(\breve{S}) = \breve{S}\times \C$, the loop $\gamma$
lifts to a loop $\gamma'\subset X \setminus L$. The loop $\gamma'$
is homotopic to zero in $X\setminus L$ hence $\gamma$ is homotopic
to zero in $S^*$. This shows that $\pi_1 (S^*) = 0$.

Now, by the Hurewicz theorem, $H_2(S^*)$ is isomorphic to the second
homotopy group of $S^*$. An element of $\pi_2 (S^*)$ can be viewed
as a continuous map $\Upsilon$ from the 2--sphere $S^2$ to $S^*$.
Without loss of generality, one may assume that its image meets
$\G$ at a finite number of general points and that
$\{\zeta_1,\ldots,\zeta_n\}= \Upsilon^{-1} (\G )$ is finite.
Consider small discs $\Delta_i$ in $S^2$ centered at $\zeta_i$.
Let $\cS_i$ be the germ of $S$ at $\Upsilon(\zeta_i)$. According
to [Ka03, Lemma 4.1], there is a germ $\cP_i \subset X$ of a
surface such that $\cS_i$ is a homeomorphic image of $\cP_i$
under $\rho$. Put
\[
\Upsilon_i = (\rho |_{\cP_i})^{-1} \circ \Upsilon |_{{\overline
\Delta_i}}\quad\text{and}\quad
S^2_0 = S^2\; \setminus\; \bigsqcup_{i=1}^n \,\Delta_i,
\]
then $\Upsilon (S^2_0)\subset \breve{S}$. By the Tietze extension
theorem, there is a continuous map $\Upsilon_0 : S^2_0 \to
\breve{X} \simeq \breve{S } \times \C$ such that $ \rho \circ
\Upsilon_0 = \Upsilon |_{S^2_0}$ and $\Upsilon_0|_{ \de \Delta_i}
=\Upsilon_i |_{\de \Delta_i}$ for every $i=1, \ldots , n$. Hence
$\Upsilon_0$ and $\Upsilon_i$'s together define a continuous map
$\Upsilon': S^2 \to X$ such that $\rho \circ \Upsilon' =
\Upsilon$. As $\pi_2 (X\setminus L) = 0$ we see that $\pi_2 (S^*)$
and hence $H_2(S^*)$ are trivial.
\end{proof}

Let $s_1, \ldots, s_k$ be the singular points of $S$. For each $i
= 1,\ldots, k$, there exists a neighborhood $U_i$ of $s_i$ in $S$
such that $U_i$ is an open cone over a closed oriented 3-manifold
$\Sigma_i = \de U_i$.  If $S\hookrightarrow\C^n$ is a closed
embedding, one can find a closed ball $B \subset \C^n$ of
sufficiently large radius such that, if $U_0 = S\setminus B$ then
$S \setminus U_0$ is a deformation retract of $S$. Put $U =
\bigsqcup_{i=0}^k\, U_i$. Note that $S_0:= S\setminus U$ is a
deformation retract of $S^*$; in particular, $\pi_1 (S_0) = H_2
(S_0) = 0$. Let $\Sigma_0 = \de U_0$ and $\Sigma = \de S_0$ so
that $\Sigma = \bigsqcup_{i=0}^k\, \Sigma_i$.

\begin{lemma}\label{lemma2.3}
Let $\Sigma$ be as above. Then $H_1(\Sigma) = H_2(\Sigma) = 0$, that
is, each of the $\Sigma_0, \ldots, \Sigma_k$ is a homology sphere.
Moreover, the 3-cycles $\Sigma_1,\ldots, \Sigma_k$ form a free basis
of $H_3(S_0) = \Z^k$.
\end{lemma}

\begin{proof}
Since $H_1 (S_0) = H_2 (S_0) = 0$ by Lemma \ref{lemma2.2}, the
exact homology sequence of pair
\[
\ldots \to H_3(S_0,\Sigma) \to H_2(\Sigma) \to H_2(S_0) \to
H_2(S_0,\Sigma) \to H_1(\Sigma) \to H_1(S_0)
\]
implies that $H_1 (\Sigma) = H_2 (S_0,\Sigma)$. By Lefschetz
duality, $H_2 (S_0,\Sigma) = H^2 (S_0)$. The latter group
vanishes because $H^2 (S_0) = \operatorname{Hom}\,(H_2\,(S_0),
\Z) = 0$, see Lemma \ref{lemma2.2}. By Poincar\'e duality,
$H_2 (\Sigma) = H_1 (\Sigma) = 0$. As $H_3 (S_0, \Sigma) =
H^1 (S_0) = 0$ and $H_4 (S_0, \Sigma) = H^0 (S_0) =\Z$,
extending the homology sequence to the left we get $0\to \Z
\to H_3(\Sigma )\to H_3(S) \to 0$. This yields the last claim.
\end{proof}

\begin{lemma}\label{lemma2.4}
Let $\cS$ be the germ of a normal surface at a point $s$, and $\cP$
the germ of a smooth surface at a point $p$. Let $\psi: \cP\to \cS$
be a finite morphism such that $\psi^{-1} (s) = p$.  Then $s$ is at
worst a quotient singularity.
\end{lemma}

\begin{proof}
Suppose that $\cS$ is embedded into $\C^n$ so that $s$ coincides
with the origin. Let $d$ be the distance function from $s$ on
$\C^n$. Put $\cS_{\varepsilon} = d^{-1}([0, \varepsilon ))\cap
\cS$, $\cS^*_{\varepsilon} =\cS_{\varepsilon} \setminus s$, $D =
d \circ \psi$, $\cP_{\varepsilon} = D^{-1}([0, \varepsilon ))$,
and $\cP_{\varepsilon}^* =\cP_{\varepsilon} \setminus p$.
According to [Mil68], $D$ has finitely many critical values.
Thus by [Mil63], there exists $\varepsilon$ such that for each
$0 < \varepsilon' \leq \varepsilon$, $\cP_{\varepsilon'}$ is a
deformation retract of $\cP_{\varepsilon}$ (resp.
$\cS_{\varepsilon'}$ is a deformation retract of
$\cS_{\varepsilon}$). We will call $\cP_{\varepsilon}$ (resp.
$\cS_{\varepsilon}$) a ``good'' neighborhood of $p$ (resp. $s$).

The fundamental group of $\cP_{\varepsilon}^*$ (resp.
$\cS_{\varepsilon}^*)$ is called the local fundamental group of
$\cP$ at $p$ (resp. $\cS$ at $s$). It is independent of the choice
of a good neighborhood, hence $\pi_1 (\cP_{\varepsilon}^*) = 0$ as
$p$ is a smooth point. Since $\psi$ is finite, there are curves
$\G_{\varepsilon} \subset \cS_{\varepsilon}$ and $Z_{\varepsilon}
\subset \cP_{\varepsilon}$ such that $ \psi |_{\cP_{\varepsilon}
\setminus Z_{\varepsilon}}: \cP_{\varepsilon} \setminus
Z_{\varepsilon} \to \cS_{\varepsilon} \setminus \G_{\varepsilon}$
is a covering, i.e. $G = \psi_*(\pi_1 (\cP_{\varepsilon} \setminus
Z_{\varepsilon}))$ is a subgroup of finite index in $\pi_1 (
\cS_{\varepsilon} \setminus \G_{\varepsilon} )$. Let $\gamma_1,
\ldots , \gamma_n$ be loops in $\cS_{\varepsilon}\setminus
\G_{\varepsilon}$ such that the cosets $[\gamma_1 ]G , \ldots ,
[\gamma_n]G$ exhaust $\pi_1 (\cS_{\varepsilon} \setminus
\G_{\varepsilon})$. Let $\gamma$ be a loop in
$\cS_{\varepsilon}^*$. After a small homotopy if necessary, we may
assume that $\gamma$ is disjoint from $\G_{\varepsilon}$. Then
$[\gamma \gamma_i] \in G$ for some $i$, which implies that
$\gamma\gamma_i$ represents zero element of $\pi_1
(\cS_{\varepsilon}^*)$ as $\pi_1 (\cP_{\varepsilon}^*) = 0$. Thus
$\pi_1 (\cS_{\varepsilon}^*)$ is finite. According to [Pr67,
Br67], $s$ may be at worst a quotient singularity.
\end{proof}

\begin{proposition}\label{prop2.5} $\empty$

\begin{enumerate}
\item[(1)] For every nontrivial algebraic $\C_+$--action on a
smooth contractible affine algebraic threefold $X$ the quotient
$S=X\|\C_+$ has at worst quotient singularities of type
$x^2+y^3+z^5=0$.
\item[(2)] $S$ is contractible.
\item[(3)] If the Kodaira logarithmic dimension $ \bar{\kappa}
(S^*)$ of $S^*$ is 1 then $S$ is smooth, and if $ \bar{\kappa}
(S^*) =-\infty$ then $S \simeq \C^2$.
\end{enumerate}
\end{proposition}

\begin{proof}
We know from Lemma \ref{lemma2.1} that $\rho: X \to S$ is
surjective and that the fibers of $\rho$ are curves.  Therefore,
we can choose a germ $\cP$ of a smooth surface at a smooth point
$p$ of $\rho^{-1}(s)$ (where $s \in S$) transversal to the curve
$\rho^{-1}(s)$. The restriction of $\rho$ to $\cP$ yields a finite
morphism $\psi: \cP \to \cS$ where $\cS$ is the germ of $S$ at
$s$. By Lemma \ref{lemma2.4}, $s$ is at most a quotient
singularity; in particular, its local fundamental group is finite.
On the other hand, by Lemma \ref{lemma2.3}, the local first
homology group at $s$ is trivial. Therefore, the local fundamental
group is perfect. The only quotient singularity whose fundamental
group is perfect is $E_8$, i.e. it is of the type $x^2 + y^3 + z^5
= 0$, see [Br67].

To prove the second statement note that $\pi_1 (S) = 0$ because
$\pi_1 (S^*) = 0$ by Lemma \ref{lemma2.2}. The statement will
follow from the Whitehead and Hurewicz theorems as soon as we show
that $H_2(S) = 0$ (since we already know that $H_i (S) = 0$ for $i
\geq 3$, see [Na67]). Let $U, \Sigma$, and $S_0$ be as defined
right before Lemma \ref{lemma2.3} so that $S = \bar{U} \cup S_0$
and $\Sigma = \bar{U}\cap S_0$. Recall that each component of $U$
is contractible, in particular, $H_2 (U) = 0$. Then $H_2(S_0) = 0$
by Lemma \ref{lemma2.2}, and $H_1(\Sigma ) = 0$ by Lemma
\ref{lemma2.3}. The Mayer--Vietoris sequence now implies that $H_2
(S) = 0$.

If $\bar{\kappa}(S^*) = 1$, any singularity of $S$ must be cyclic
quotient [GuMiy92], hence $S$ is smooth by virtue of (1). If
$\bar{\kappa} (S^*) = -\infty$, the only logarithmic contractible
surfaces with at worst $E_8$--type singularities are $\C^2$ or the
surface $x^2+y^3+z^5 = 0$ in $\C^3$, see [MiySu91, Theorem 2.7].
The second possibility should be eliminated because $\pi_1 (S^*)
\ne 0$ contrary to Lemma \ref{lemma2.2}. This implies (3).
\end{proof}

\begin{corollary}\label{cor2.6}
Every smooth contractible affine algebraic threefold with a
nontrivial algebraic $\C_+$-action is rational.
\end{corollary}

\begin{proof}
According to Proposition \ref{prop2.5}, surface $S$ is contractible
logarithmic (i.e. it has at worst quotient singularities) and hence
rational by [GuPrSh97, PrSh97, GuPr99]. Therefore, $S \times \C$ is
rational, and so is $X$ since it is an affine modification of $S
\times \C$, see [Ka03, Lemma 2.2].
\end{proof}

\begin{theorem}\label{thm2.7}
For every nontrivial algebraic $\C_+$--action on a
smooth contractible affine algebraic threefold $X$, the quotient
$S=X\|\C_+$ is a smooth contractible affine surface.
\end{theorem}

\begin{proof}
Let $S_0, \Sigma$, and $\Sigma_i$ be as defined right before Lemma
\ref{lemma2.3}. Assume first that $S$ has only one singular point.
Then the boundary $\Sigma$ of $S_0$ consists of two components.
One component is $\Sigma_1$, which is the link of singularity at 0
of $x^2 + y^3 + z^5 = 0$, according to Proposition \ref{prop2.5}.
The manifold $\Sigma_1$ is also known as the Poincar\'e homology
sphere. The other component is $\Sigma_0$, which is also a homology
sphere by Lemma \ref{lemma2.3}. Lemmas \ref{lemma2.2} and
\ref{lemma2.3} imply that $\pi_1 (S_0) =0$ and that the embeddings
$\Sigma_0\hookrightarrow S_0$ and $\Sigma_1\hookrightarrow S_0$
induce isomorphisms in homology. Thus $S_0$ is a simply connected
homology cobordism between $\Sigma_1$ and $\Sigma_0$. But this
contradicts the Fintushel--Stern theorem [FiSt90, Theorem 5.2]
which says that the Poincar\'e homology sphere cannot be homology
cobordant to any other homology sphere via a simply connected
homology cobordism.


To complete the proof, it is enough to consider the case of two
singular points; the general case will follow by a similar
argument. If $S$ has two singular points, $\Sigma$ is a disjoint
union of $\Sigma_0$, $\Sigma_1$, and $\Sigma_2$. Let us join a
point $x_0\in \Sigma_0$ with a point $x_2 \in \Sigma_2$ by a path
$\gamma$ in $S_0$. Let $V_2$ and $V_1$ be tubular neighborhoods of
$\gamma$ in $S_0$ (i.e. each $V_i$ is homeomorphic to $\gamma
\times B_i$ where $B_i$ is a three-dimensional ball, and $V_i$
meets $\Sigma_j$, $j=0, 2$, along the ball $x_j \times B_i$) such
that ${\rm int} \, V_2 \supset V_1$. Put $S_1 = S_0\setminus V_1$.
Then the boundary of $S_1$ consists of two components, $\Sigma_1$
and $\Sigma'$, where $\Sigma'$ is a connected sum of $\Sigma_0$
and $\Sigma_2$ (and hence is a homology sphere). Note that $\pi_1
(S_1) = \pi_1 (S_0\setminus \gamma) = 0$ by the dimension
argument. In order to show that we have a homology cobordism
between $\Sigma_1$ and $\Sigma'$ and thus get a contradiction with
the Fintushel--Stern theorem, we only need to show that $H_2 (S_1)
= 0$ and the 3-cycle $\Sigma_1$ generates $H_3(S_1)=\Z$.

The Mayer--Vietoris sequence of $S_0 =V_2\,\cup\,S_1$ implies
that $H_2 (S_1)$ is the image of $H_2 (V_2\setminus V_1)$ under
the natural embedding. Note that $x_2 \times (B_2 \setminus B_1)$
is a deformation retract of $V_2 \setminus V_1$. Therefore, every
element of $H_2 (S^1)$ can be represented by a 2-cycle in $x_2
\times (B_2 \setminus B_1) \subset \Sigma_2 \setminus (x_2 \times
B_1)$. As $\Sigma_2$ is a homology sphere, we conclude that $H_2
(\Sigma_2 \setminus B_2 ) = 0$ and hence $H_2 (S_1) = 0$. As
$H_3(V_2\setminus V_1)=H_3(V_2)=0$ and $H_2(V_2\setminus V_1)=\Z$,
applying the Mayer--Vietoris sequence again we get exact sequence
$0\to H_3(S_1)\to H_3(S_0)\to \Z \to 0$. Since $\{ \Sigma_1 ,
\Sigma_2 \}$ is a free basis of $H_3(S_0)$ according to Lemma
\ref{lemma2.3}, we see that $H_3(S_1)$ is freely generated by
$\Sigma_1$.

This leaves us with just one possibility that $S$ has no singular
points and hence is smooth. That it is contractible was already
proved in Proposition \ref{prop2.5} (2).
\end{proof}

\begin{corollary}\label{cor2.8}
Let $X$ be a smooth contractible affine algebraic threefold with a
nontrivial algebraic $\C_+$-action on it.
\begin{enumerate}
\item[(1)] If the action is free then $X$ is isomorphic to $S
\times \C$ and the action is induced by a translation on the
second factor.
\item[(2)] If $X$ admits a dominant morphism from a threefold of
form $C \times \C^2$ then the algebraic quotient $S=X\|\C_+$ is
isomorphic to $\C^2$.
\item[(3)] If the assumptions of both (1) and (2) hold then $X$ is
isomorphic to $\C^3$.
\end{enumerate}
\end{corollary}

\begin{proof}
The first statement was proved in [Ka03, Theorem 5.4 (ii)] under
the additional assumption that the $S = X\|\C_+$ is smooth.
Theorem \ref{thm2.7} removes this assumption and proves (1) in
full generality. In the second statement, we have a dominant
morphism $C \times \C \to S$. As the Kodaira logarithmic dimension
$\bar{\kappa} (C \times \C )$ equals $-\infty$ we conclude that
$\bar{\kappa} (S) = -\infty$. Since $S$ is also smooth and
contractible, it is isomorphic to $\C^2$ (see e. g. [Miy01]). The
third statement is an obvious consequence of (1) and (2).
\end{proof}

Two $\C_+$-actions on a variety are said to be equivalent if
they have the same general orbits (or, equivalently, the
corresponding locally nilpotent derivations have the same kernel).
In particular, non-equivalent actions generate different quotient
morphisms. Corollary \ref{cor2.8} (3) implies the following result.

\begin{corollary}\label{cor2.9}
Suppose that a smooth contractible affine algebraic threefold $X$
admits two non-equivalent nontrivial algebraic $\C_+$--actions.
Then $X\|\C_+=\C^2$ for any nontrivial algebraic $\C_+$-action.
Furthermore, $X$ is isomorphic to $\C^3$ if it admits a free
$\C_+$-action.
\end{corollary}

It is worth mentioning that Corollary \ref{cor2.6} also follows from
Theorem \ref{thm2.7} and [GuSh89].

\section{The case when $S\simeq \C^2$}

The aim of this section is to give an independent proof of Corollary
\ref{cor2.8} (2) (and hence of Corollary \ref{cor2.8} (3)) which does
not use the Fintushel-Stern theorem.

Let $X$ be the complement to an effective divisor $D$
of simple normal crossing type in a projective algebraic manifold
$\bX$.
Consider the sheaf $\Omega^k (\bar {X}, D)$ of logarithmic
$k$-forms on $\bX$ along $D$ (that is, each section of this sheaf
over an open subset $U\subset \bX$ is a holomorphic $k$-form on
$U\cap X$ which has at most simple poles at general points of $U
\cap D$). Let $r$ be the rank of $\Omega^k (\bar {X}, D)$,
$S^m\Omega^k (\bar {X}, D)$ its symmetric $m$-power, and $H_0
(\bX, S^m\Omega^k (\bar {X}, D))$ the space of holomorphic
sections of $\Omega^k (\bar {X}, D)$ over $\bX$. We say that the
Kodaira--Iitaka--Sakai logarithmic $k$-dimension $
\bar{\kappa}_k(X)$ of $X$ is $ - \infty$ if no symmetric power of
$\Omega^k (\bar {X}, D)$ has a non-trivial global section, and
otherwise we put
\[
\bar{\kappa}_k(X)= \limsup_{m \to + \infty } {\frac {\log \dim
H_0(\bX , S^m\Omega^k (\bar {X}, D))}{\log m}} -r +1.
\]
This definition does not depend on the choice of simple normal
crossing completion $\bX$ of $X$, see [Ii77, Ka99]. One can
easily see that $ \bar{\kappa}_k(X)=-\infty$ if $k >\dim X$, and
$\bar{\kappa}_k(X)$ is the usual Kodaira logarithmic dimension
in the case when $k =\dim X$. \vs

\begin{lemma}\label{lemma3.1}{\rm [Ii77, Ka99, Prop. 4.2]}
Let $\bX_1$ and $\bX_2$ be complete complex algebraic manifolds,
and $D_1$ and $D_2$ divisors of SNC-type in $\bX_1$ and $\bX_2$,
respectively. Suppose that $\bar{f}: \bX_1 \to \bX_2$ is a
morphism and that $\bar{f}$ is an extension of a dominant
morphism $f: X_1 \to  X_2$ where $X_i = \bX_i - D_i$. Then
$\bar{f}$ generates a natural homomorphism $f^*: S^m \Omega^k
(\bX_2, D_2) \to S^m \Omega^k (\bX_1 , D_1)$.
\end{lemma}

The word ``natural" above means that we treat $H_0(\bX_i,S^m
\Omega^k(\bX_i, D_i))$ as the subspace of $H_0(\bX,\Omega^k(\bX_i,
D_i)^{\otimes m})$ invariant under the natural action of the
symmetric group $S(m)$ and that $f^*$ is generated by the induced
mapping of $k$--forms. In particular, $f^*$ sends nonzero sections
of $S^m\Omega^k (\bX_2, D_2)$ to nonzero sections of $S^m \Omega^k
(\bX_1 , D_1)$. Therefore we have the following result.

\begin{corollary}\label{cor3.2}
Let $f: X_1 \to X_2$ be a dominant morphism of algebraic varieties
and $n_i =\dim X_i$. Then $\bar{\kappa}_k (X_1) + C_{n_1,k} \geq
\bar{\kappa}_k (X_2)+C_{n_2,k}$ where $C_{n_i,k}$ is the number of
combinations. In particular, if $\bar{\kappa}_k (X_1) = -\infty$
then $\bar{\kappa}_k(X_2) = -\infty$.
\end{corollary}

Let $H$ be a hyperplane in $\P^s$, i.e. $\C^s=\P^s \setminus H$.
Then $X'=\bX \times \P^s$ is a completion of $X \times \C^s$ and
$D'=X' \setminus (X \times \C^s)$ is of simple normal crossing
type. Using the fact that any sheaf of form
$$\Omega^{1}(\P^s,H)^{\otimes m_1} \otimes \cdots \otimes
\Omega^{s}(\P^s, H)^{\otimes m_s}$$ has no global nonzero sections
over $\P^s$, one can show that $$H_0(\bX ,\Omega^k(\bX
,D)^{\otimes m})= H_0(X',\Omega^k(X',D')^{\otimes m})$$
which implies the following.

\begin{lemma}\label{lemma3.3}
Let $X'=X \times \C^s$ and $n = \dim X$. Then $\bar{\kappa}_k (X')
= \bar{\kappa}_k (X)+C_{n,k}-C_{n+s,k}$ for any $k\geq 0$. In
particular, $\bar{\kappa}_k (X') = -\infty$ when $k
>n$.
\end{lemma}

Applying the theorem about removing singularities of holomorphic
functions in codimension 2 we get the following result.

\begin{lemma}\label{lemma3.4}
Let $Z$ be a subvariety of codimension at least 2
in a algebraic manifold $X$. Then $ \bar{\kappa}_k (X)
=\bar{\kappa}_k (X \setminus Z)$ for every $k$.
\end{lemma}

\begin{theorem}\label{thm3.5}
Let $X$ be a smooth contractible affine algebraic threefold such
that $\bar{\kappa}_2(X) = -\infty$. Then, for every nontrivial
algebraic $\C_+$-action on $X$, the algebraic quotient $S =
X\|\C_+$ is isomorphic to $\C^2$.
\end{theorem}

\begin{proof}
Let $F$ be the set of singular points of $S$. According to Lemma
\ref{lemma2.1}, $L= \rho^{-1}(F)$ is a curve. Therefore,
$\bar{\kappa}_2 (X\setminus L) = -\infty$, see Lemma
\ref{lemma3.4}. By Corollary \ref{cor3.2}, $\bar{\kappa}_2 (S^*) =
-\infty$, and the statement follows from Proposition \ref{prop2.5}
(3).
\end{proof}

Now Lemma \ref{lemma3.3} implies Corollary \ref{cor2.8} (2).

\begin{remark}\label{rem3.6}
Consider an $n$-dimensional smooth contractible affine algebraic
variety $X$ and a free algebraic action of a unipotent group $U$
on it (that is, each orbit of the action is isomorphic to $U$).
Suppose that $U$ is of dimension $n-2$ (i.e. $U$ is isomorphic to
$\C^{n-2}$ as an affine algebraic variety). It was mentioned in
[Ka03, Remark 5.5] that the morphism $X \to S = X\|U$ is
surjective. Since surjectivity of the quotient morphism is the
only crucial argument in the proof of Proposition 2.5, we can
extend our results to the action of $U$. That is, if $X$ admits a
dominant morphism from an $n$--fold of form $C \times \C^{n-1}$
then $S$ is isomorphic to $\C^2$, $X$ is isomorphic to $\C^n
\simeq S \times U$, and the action of $U$ is induced by the
natural action on the second factor.
\end{remark}

\end{document}